\documentclass[11pt,a4paper]{article}
\usepackage[numbers,sort&compress]{natbib}
\usepackage[colorlinks,
            linkcolor=blue,
            anchorcolor=green,
            citecolor=magenta
            ]{hyperref}
\usepackage{mathrsfs,amssymb,amsfonts}
\usepackage{graphicx}
\usepackage{amsmath,amssymb,latexsym,color}
\usepackage[mathscr]{eucal}
\usepackage[numbers]{natbib}
\usepackage{bm}
\usepackage[top=1in, bottom=1in, left=1.25in, right=1.25in]{geometry}
\makeatletter

\newcommand{\Rmnum}[1]{\expandafter\@slowromancap\romannumeral #1@}
\makeatother
\newtheorem{thm}{Theorem}[section]
\newtheorem{defin}[thm]{Definition}

\newtheorem{lem}[thm]{Lemma}

\newtheorem{example}{Example}
\newtheorem{remark}{Remark}[section]

\newcommand{\qed}{\hfill\Box\medskip}
\usepackage{CJK}
\usepackage{geometry}
\begin{document}
\begin{CJK*}{GBK}{song}
 \setlength{\baselineskip}{14pt}
\renewcommand{\abovewithdelims}[2]{
\genfrac{[}{]}{0pt}{}{#1}{#2}}

\title{\bf Edge-fault-tolerant edge-bipancyclicity of balanced hypercubes }

\author{Pingshan Li \quad  Min Xu\footnote{Corresponding author. \newline {\em E-mail address:} xum@bnu.edu.cn (M. Xu) .}\\
{\footnotesize   \em  Sch. Math. Sci. {\rm \&} Lab. Math. Com. Sys., Beijing Normal University, Beijing, 100875,  China} }
 \date{}
 \date{}
 \maketitle

\begin{abstract}
The balanced hypercube, $BH_n$, is a variant of hypercube $Q_n$. Hao et al. [Appl. Math. Comput. 244 (2014) 447-456] showed that there exists a fault-free Hamiltonian path between any two adjacent vertices in $BH_n$ with $(2n-2)$ faulty edges.
Cheng et al. [Inform. Sci. 297 (2015) 140-153] proved that $BH_n$ is $6$-edge-bipancyclic  after  $(2n-3)$ faulty edges occur for all $n\ge2$.
In this paper,  we improve these two results by  demonstrating that $BH_n$ is $6$-edge-bipancyclic even when there exist $(2n-2)$ faulty edges for all $n\ge2$. Our result is optimal with respect to the maximum number of tolerated edge faults.


\medskip
\end{abstract}

\section{Introduction}
In the field of parallel and distributed systems, interconnection networks are an important research area. Typically, the topology
of a network can be represented as a graph in which the vertices represent processors and the edges represent communication links.

The hypercube network has been proved to be one of the most popular interconnection networks as it possesses many excellent properties such as a recursive structure, regularity, and symmetry. It is well known that no network typically meets all the aspects of a given set of requirements. Thus, a number of hypercube variants have been proposed, such as folded hypercubes \cite{A. Elamawy}, crossed cubes \cite{K. Efe}, M\"{o}bius cubes \cite{P.Cull}, twisted cubes \cite{PAJHibers}, and shuffle cubes \cite{Li} and so on(see \cite{Xujunming}).

  The balanced hypercube, proposed by Huang and Wu \cite{Huang}, is also a  hypercube variant. Similar to hypercubes, balanced hypercubes are bipartite graphs \cite{Huang} that are vertex-transitive \cite{Huang2} and edge-transitive \cite{Zhou}. Balanced hypercubes are superior to hypercubes in that they have a smaller diameter as compared to hypercubes.

   Studies on balanced hypercubes can be found in
\cite{Dongqincheng, Dongqincheng2, R.X.Hao, Huang, Huang2, H.Lu, Huazhonglu, Xumin, Mingchengyang, Mingchengyang2, Qingguozhou, Zhou}.

For graph definitions and notations, we follow \cite{J.A.Bondy}.
A graph $G$ consists of a vertex set $V(G)$ and an edge set $E(G)$, where an edge is an unordered pair of distinct vertices of $G$. A graph $G$ is called bipartite if its vertex set can  be partitioned into two parts $V_1, V_2$  such that every edge has one endpoint in $V_1$ and one in $V_2$. A vertex $v$ is a neighbor of $u$ if $(u, v)$ is an edge of $G$,
and $N_{G}(u)$ denotes the set of all the neighbors of $u$ in $G$.
A path  $P$ of length $\ell$ from $x$ to $y$, denoted by $\ell$-path $P$, is a finite sequence of distinct vertices $\langle v_0, v_1, \cdots, v_{\ell}\rangle$ such that $x=v_0, y=v_{\ell}$, and $(v_i, v_{i+1})\in E$ for $0\le i\le \ell-1$.
We also denote the path $P$ as $\langle v_0, v_1, \cdots, v_i, Q, v_j, v_{j+1}, \cdots, v_{\ell}\rangle$, where $Q$ is the path $\langle v_i, v_{i+1}, \cdots, v_j\rangle$.
A cycle $C$ of length $\ell+1$ is a closed path$\langle v_0, v_1, \cdots, v_{\ell}, v_0\rangle$, denoted by $(\ell+1)$-cycle $C$.

In an interconnection network, the problem of simulating one network by another is modeled as a graph embedding problem. In all embedding problems, the cycle embedding problem is one of the most common problem; it refers to finding a cycle of a given length in a graph. A graph $G$ of order $|V(G)|$ is $m$-pancyclic,
if it contains  every $\ell$-cycle for $m\le \ell\le |V(G)|$. A bipartite graph $G$ is $m$-bipancyclic,
 if it contains every even $\ell$-cycle for $m\le \ell\le |V(G)|$. A graph $G$ is pancyclic (resp. bipancyclic) if it is $g$-pancyclic ($g$-bipancyclic), where $g=g(G)$ is the girth of $G$.
 A graph $G$ is vertex-pancyclic (resp. edge-pancylic) if every vertex (resp. edge) lies on various $\ell$-cycles for all $g\le \ell\le V(G)$. A path is called a Hamiltonian path if it contains all the vertices of $G$. A graph $G$ is said to be Hamiltonian connected if there exists a Hamiltonian path between any two vertices of $G$. A bipartite graph is Hamiltonian laceable if there is a Hamiltonian path between any two vertices in different bipartite sets.

A bipartite  graph $G$ is $k$-fault-tolerant hamiltonian laceable (resp. bipancyclic, vertex-bipancyclic, and edge-bipancyclic ) if $G-F$ remains Hamiltonian laceable  (resp. bipancyclic, vertex-bipancyclic, and edge-bipancyclic ) for $F\subseteq V(G)\cup E(G)$, $|F|\le k$.
A  bipartite graph $G$ is $k$-edge-fault-tolerant Hamiltonian laceable  (resp. bipancyclic, vertex-bipancyclic, and edge-bipancyclic ) if $G-F$ remains Hamiltonian $lacelabe$ (resp. bipancyclic, vertex-bipancyclic, and edge-bipancyclic ) for $F\subseteq E(G), |F|\le k$.

The balanced hypercube, $BH_n$, has been studied by many researchers.
 Xu  et al. \cite{Xumin} proved that $BH_n$ is edge-bipancyclic and Hamiltonian laceable.
 Yang \cite{Mingchengyang} proved that  $BH_n$ is bipanconnected.
 Yang  \cite{Mingchengyang2} also demonstrated that the super connectivity of $BH_n$ is $(4n-4)$ and the super edge-connectivity of $BH_n$ is $(4n-2)$ for $n\ge 2$.
 L\"{u} et al. \cite{Huazhonglu} proved that $BH_n$ is hyper-Hamiltonian laceable.
 Cheng et al. \cite{Dongqincheng} proved that $BH_n$ is $(n-1)$-vertex-fault-tolerant edge-bipancyclic.
  Hao et al. \cite{R.X.Hao} showed that there exists a fault-free Hamiltonian path between any two adjacent vertices in $BH_n$ with $(2n-2)$ faulty edges.
  Zhou et al. \cite{Qingguozhou} proved that $BH_n$ is $(2n-2)$-edge-fault-tolerant Hamiltonian laceable.
 Cheng et al. \cite{Dongqincheng2}  proved that $BH_n$ is $(2n-3)$ edge-fault-tolerant $6$-edge-bipancyclic  for all $n\ge2$.
 In this paper, we improve the results of  Hao et al. \cite{R.X.Hao} and  Cheng et al. \cite{Dongqincheng2} by demonstrating that $BH_n$ is $(2n-2)$ edge-fault-tolerant $6$-edge-bipancyclic  for all $n\ge2$. Our result is optimal with respect to the maximum number of tolerated edge faults.

The rest of this paper is organized as follows. In Section 2, we introduce  two equivalent definitions  of balanced hypercubes  and discuss some of their properties. In Section 3, we investigate edge-bipancyclic of $BH_n$ with faulty edges. Finally, we conclude this paper in Section $4$.

\section{Balanced hypercubes}

\begin{figure}[!htbp]
  \centering
  \includegraphics[width=0.6\textwidth]{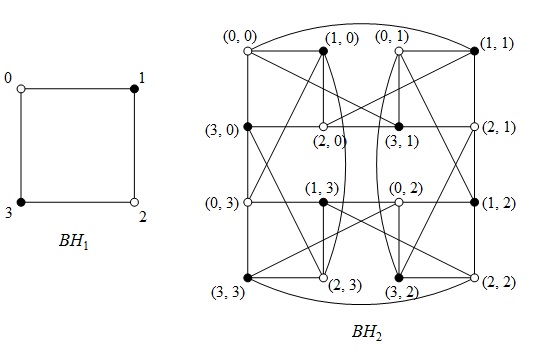}\\
  \caption{Illustration of $BH_1$ and $BH_2$}\label{BH2}
\end{figure}

Wu and Huang  \cite{Huang} presented two equivalent definitions of $BH_n$ as follows:

\begin{defin}  An $n$-dimensional balanced hypercube $BH_n$ has $2^{2n}$ vertices, each labeled by an $n$-bit string  $(a_0, a_1, \cdots, a_{n-1})$, where $a_i\in \{0, 1, 2, 3\}$ for all $0\le i\le n-1$. A arbitrary vertex $(a_0, a_1, \cdots, a_{i-1}, a_{i}$, $a_{i+1}, \cdots, a_{n-1})$ is adjacent to the following $2n$ vertices:
\begin{center}
  $\begin{array}{l}
          (1)~~ ((a_0\pm1)\mod\ 4, a_1, \cdots, a_{i-1}, a_{i}, a_{i+1}, \cdots, a_{n-1})~where~1\le i\le n-1, \vspace{0.5ex} \\
           (2)~~ ((a_0\pm1)\mod\ 4, a_1, \cdots,  a_{i-1}, (a_{i}+(-1)^{a_0})\mod\ 4, a_{i+1}, \cdots, a_{n-1})~where~1\le i\le n-1.
         \end{array}
  $
 \end{center}
In $BH_n$, the first coordinate $a_0$ of vertex $(a_0, a_1, \cdots, a_{n-1})$ is called the inner index, and the second coordinate $a_i(1\le i\le n-1)$ is called the $i$-dimension index. From the definition, we have that $N_{BH_n}((a_0, a_1, \cdots, a_{n-1}))=N_{BH_n}((a_0+2, a_1, \cdots, a_{n-1}))$.
Figure \ref{BH2} shows two balanced hypercubes of dimensional one and two.
\end{defin}

 Briefly, we assume that `$+, -$' for the coordinate of a vertex is an operation with mod 4 in the remainder of the paper.
 Let $X_{j, i}=\{(a_0, a_1, \cdots, a_{j-1}, a_{j}, a_{j+1}, \cdots, a_{n-1})\mid a_k\in \{0, 1, 2, 3\}, 0\le k\le n-1, a_j=i\}$ for $1\le j\le n-1$ and $ i\in \{0, 1, 2, 3\}$ and let $BH^{j, i}_{n-1}=BH_n[X_{j, i}]$. Then, $BH_n$ can be divided into four copies: $BH^{j, 0}_{n-1}, BH^{j, 1}_{n-1}, BH^{j, 2}_{n-1}, BH^{j, 3}_{n-1}$ where $BH^{j, i}_{n-1}\cong BH_{n-1}$ for $i=0, 1, 2, 3$ \cite{Dongqincheng}. We use $BH^i_{n-1}$ to denote  $BH^{n-1, i}_{n-1}$ for $i=0,1, 2, 3$ .

\begin{defin} The balanced hypercube $BH_n$ can be constructed recursively as follows:
\begin{enumerate}
  \item $BH_1$ is a $4$-cycle with vertex-set $\{0, 1, 2, 3\}$.
  \item $BH_{n}$ is a construct from four copies of $BH_{n-1}: BH^0_{n-1}, BH^1_{n-1}, BH^2_{n-1}, BH^3_{n-1}$. Each vertex $(a_0, a_1, \cdots$, $ a_{n-2}, i)$ has two extra adjacent vertices:

$\begin{array}{l}
  ~(1)~ In~ BH^{i+1}_{n-1}: (a_0 \pm 1, a_1, \cdots,  a_{n-2}, i+1)~if~a_0~is~ even\vspace{1.0ex}. \\
   ~(2)~ In~ BH^{i-1}_{n-1}: (a_0 \pm 1, a_1, \cdots, a_{n-2}, i-1)~if~a_0~is~ odd.
 \end{array}
$

\end{enumerate}

\end{defin}

Since $BH_n$ is a bipartite graph, then $V(BH_n)$ can be divided into two disjoint parts. Obviously, the vertex-set $V_1=\{a=(a_0, a_1, \cdots, a_{n-1})\mid a\in V(BH_n)$ and $ a_0$  is odd$\}$ and $V_2=\{a=(a_0, a_1, \cdots, a_{n-1})\mid a\in V(BH_n)$ and $ a_0$  is even$\}$ form the desired partition. We use black nodes to denote the vertices in $V_1$ and white nodes to denote the vertices in $V_2$.

Let $(u, v)$ be an edge of $BH_n$, if $u$ and $v$ differ only with regard to the inner index, then $(u, v)$ is  said to be a 0-dimensional edge. If $u$ and $v$ differ not only in terms of the inner index but also with regard to the $i$-dimension index, then $(u, v)$ is called the $i$-dimensional edge. We use $\partial D_d(0\le d\le n-1)$ to denote the set of all $d$-dimensional edges.

There are some known properties about $BH_n$.

\begin{lem}{\rm (\cite{Huang2, Zhou})}\label{L1} The balanced hypercube $BH_n$ is vertex-transitive and edge-transitive.
\end{lem}

\begin{lem}{\rm (\cite{Qingguozhou})}\label{LL2} The balanced hypercube $BH_n$ is $(2n-2)$-edge-fault-tolerant Hamiltonian laceable for $n\ge 2$.
\end{lem}

\begin{lem}{\rm (\cite{Xumin})}\label{L3} The balanced hypercube $BH_n$ is edge-bipancyclic for $n\ge2$.
\end{lem}

\begin{lem} {\rm ({\cite{Dongqincheng}})}\label{L4} Let $e=(x, y)$ be an arbitrary  edge in $BH^{j, 0}_{n-1}$. Then, there exist two internal vertex-disjoint paths $\langle x, x_1, y_1, x_2, y_2, x_3, y_3, y\rangle$ and $\langle x, x'_1, y'_1, x'_2, y'_2, x'_3, y'_3, y\rangle$ in $BH_n$ such that $(x_i, y_i), (x'_i, y'_i)$$\in E(BH^{j, i}_{n-1})$ where $1\le j\le n-1$ and $i=1, 2, 3$.
\end{lem}

\begin{lem}{\rm (\cite{H.Lu})}\label{LL3} Let $n\ge 2$ be an integer. Then, $BH_n-\partial D_0$  has four components, and each component is isomorphic to $BH_{n-1}$.
\end{lem}

\noindent{\bf Remark.} The above Lemma shows that one can divide $BH_n$ into four $BH_n$s by deleting $\partial D_d$ for any $d\in\{0, 1, \cdots, n-1\}$. The four components  of $BH_n$ through the deletion of $\partial D_j$ are  $BH^{j, 0}_{n-1}$, $BH^{j, 1}_{n-1}$, $BH^{j, 2}_{n-1}$, and $BH^{j, 3}_{n-1}$ for  $1\le j\le n-1$. For convenience, we use $BH^{0, 0}_{n-1}, BH^{0, 1}_{n-1}, BH^{0, 2}_{n-1}, and BH^{0, 3}_{n-1}$ to denote the components of $BH_n-\partial D_0$ throughout this paper.

A graph $G$ is hyper-Hamiltonian laceable if it is Hamiltonian laceable and, for an arbitrary vertex $v$ in $V_{i}$ where $i\in \{0, 1\}$, there exists a Hamiltonian path in $G-v$ joining any two different vertices in $V_{1-i}$.  L\"{u} et al. obtained the following result.

\begin{lem}{\rm (\cite{Huazhonglu})}\label{LL4} The balanced hypercube $BH_n$ is hyper-Hamiltonian laceable for $n\ge 1$.
\end{lem}

In the following, we discuss some properties that are used in the proof of our main results.

\begin{lem}\label{L6} For an arbitrary vertex $u$ in $BH^{j, i}_{n-1}$ where $0\le j\le n-1$, $0\le i\le 3$. Suppose that $F\subseteq E(BH_n), |F|\le 2n-2$ and $|F\cap BH^{j, i}_{n-1}|\le 2n-3$. Then, there exists a  $2$-path $\langle u, v, w\rangle\subseteq BH_n\setminus F$ where $u, v\in BH^{j, i}_{n-1}, w\in BH_n\setminus BH^{j,i}_{n-1}$.
\end{lem}

\noindent{\bf Proof: } Without loss of generality, we can assume that $u=(0, 0, \cdots, 0)\in BH^{0}_{n-1}$.
Note that $N_{BH^{0}_{n-1}}(u)=2n-2$
and  $u$ is a white vertex, there exist $2(2n-2)$ different edges from $N_{BH^{0}_{n-1}}(u)$ to
 $BH^{3}_{n-1}$. Suppose that $|F\cap BH^{0}_{n-1}|=k, |F\cap (BH_n\setminus BH^{0}_{n-1})|=t$. We have

~~~~~~~~~~~~~~~~~~~~~~~~~~~~~~~~~~~~~~~~~~~~~~~~~~~~~~~~~~~~~$\left\{
\begin{aligned}
  &k+t\le 2n-2; \\
 & k\le 2n-3.
\end{aligned}
\right.$\vspace{1.0ex}

\begin{figure}[!htbp]
  \centering
  \includegraphics[width=0.43\textwidth]{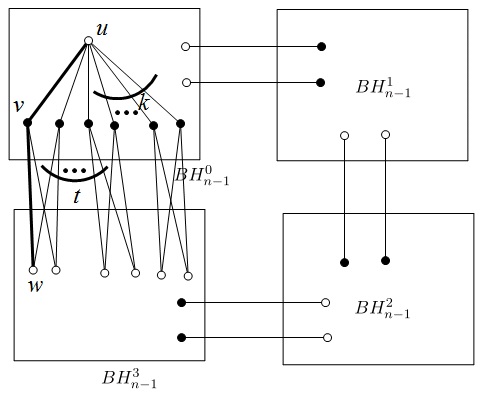}\\
  \caption{The fault-free path $\langle u, v, w\rangle$ of Lemma \ref{L6}}\label{L2_7}
\end{figure}

Hence, there exists at least one $2$-path $\langle u, v, w\rangle\subseteq BH_n\setminus F$ where $u, v\in BH^{0}_{n-1}, w\in BH^{3}_{n-1}$ owing to $2((2n-2)-k)-t\ge 2(2n-2)-(k+t)-k\ge 2n-2-k\ge 1$. See figure \ref{L2_7} for illustration.

\begin{lem}\label{L9} Suppose that $e=(u, v)$ is an edge between $BH^{j, i}_{n-1}$ and $BH^{j, i+1}_{n-1}$ where $0\le j\le n-1, 0\le i\le 3$ for $n\ge 2$. Then, there exists a cycle $C$ of length $8$ in $BH_n\setminus F$ where $F\subseteq E(BH_n)$, $|F|\le 2n-2$ and $|F\cap \partial D_{j}|\ge 1$ such that $|E(C)\cap BH^{j, i}_{n-1}|=1$.
\end{lem}

\begin{figure}[!htbp]
  \centering
  \includegraphics[width=0.4\textwidth]{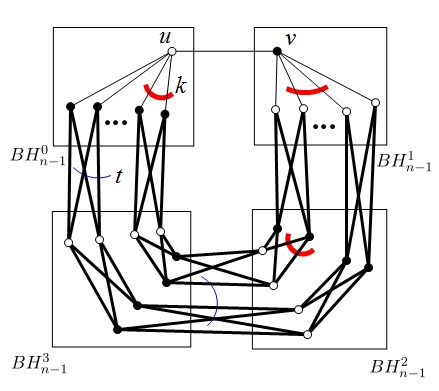}\\
  \caption{The $5$-paths  in Lemma \ref{L9}}\label{L2_10}
\end{figure}

\noindent{\bf Proof: }By Lemma \ref{L1}, $BH_n$ is edge-transitive, Without loss of generality, let $j=n-1$ and  $u=(0, 0, \cdots, 0), v=(1, 0, \cdots, 0, 1)$. There exist $4(n-1)$ edge disjoint paths of length $5$ from $N_{BH^{0}_{n-1}}(u)$ to $N_{BH^{1}_{n-1}}(v)$ such that each path has an edge in $BH^2_{n-1}$ and $BH^{3}_{n-1}$\vspace{1.0ex}. We list them as follows (see figure \ref{L2_10}):  \vspace{1.0ex}\\
$\begin{array}{l}
   P_{0, 1}=\langle  (1, 0, \cdots, 0), (2, 0, \cdots, 0, 3), (3, 0, \cdots, 0, 3), (0, 0, \cdots, 0, 2), (1, 0, \cdots, 2), (2, 0, \cdots, 1)\rangle; \\
   P_{0, 2}= \langle  (1, 0, \cdots, 0), (0, 0, \cdots, 0, 3), (1, 0, \cdots, 0, 3), (2, 0, \cdots, 0, 2), (3, 0, \cdots, 2), (0, 0, \cdots, 1)\rangle;\\
   P_{0, 3}= \langle  (3, 0, \cdots, 0), (2, 0, \cdots, 0, 3), (1, 0, \cdots, 0, 3), (0, 0, \cdots, 0, 2), (3, 0, \cdots, 2), (2, 0, \cdots, 1)\rangle;\\
   P_{0, 4}=\langle  (3, 0, \cdots, 0), (0, 0, \cdots, 0, 3), (3, 0, \cdots, 0, 3), (2, 0, \cdots, 0, 2), (1, 0, \cdots, 2), (0, 0, \cdots, 1)\rangle;\\
   P_{k, 1}=\langle  (1, \overbrace{0, \cdots, 0}^{k-1}, 1, \overbrace{0, \cdots, 0, 0}^{n-k-1}), (2, 0, \cdots, 0, 1, 0, \cdots, 0, 3), (3, 0, \cdots, 0, 2, 0, \cdots, 0, 3), \vspace{1.0ex}\\~~~~~~~~~(0, 0, \cdots, 0, 2, 0, \cdots, 0, 2), (1, 0, \cdots, 0, 3, 0, \cdots, 0, 2), (2, 0, \cdots, 0, 3, 0, \cdots, 0, 1) \rangle;  \\
   P_{k, 2}= \langle  (1, \overbrace{0, \cdots, 0}^{k-1}, 1, \overbrace{0, \cdots, 0, 0}^{n-k-1}), (0, 0, \cdots, 0, 1, 0, \cdots, 0, 3), (1, 0, \cdots, 0, 2, 0, \cdots, 0, 3), \vspace{1.0ex}\\~~~~~~~~~(2, 0, \cdots, 0, 2, 0, \cdots, 0, 2), (3, 0, \cdots, 0, 3, 0, \cdots, 0, 2), (0, 0, \cdots, 0, 3, 0, \cdots, 0, 1)\rangle;
      \end{array}$

   $\begin{array}{l}
   P_{k, 3}= \langle  (3, \overbrace{0, \cdots, 0}^{k-1}, 1, \overbrace{0, \cdots, 0, 0}^{n-k-1}), (2, 0, \cdots, 0, 1, 0, \cdots, 0, 3), (1, 0, \cdots, 0, 2, 0, \cdots, 0, 3), \vspace{1.0ex}\\~~~~~~~~~(0, 0, \cdots, 0, 2, 0, \cdots, 0, 2), (3, 0, \cdots, 0, 3, 0, \cdots, 0, 2), (2, 0, \cdots, 0, 3, 0, \cdots, 0, 1)\rangle;\\
   P_{k, 4}= \langle (3, \overbrace{0, \cdots, 0}^{k-1}, 1, \overbrace{0, \cdots, 0, 0}^{n-k-1}), (0, 0, \cdots, 0, 1, 0, \cdots, 0, 3), (3, 0, \cdots, 0, 2, 0, \cdots, 0, 3), \vspace{1.0ex}\\~~~~~~~~~(2, 0, \cdots, 0, 2, 0, \cdots, 0, 2), (1, 0, \cdots, 0, 3, 0, \cdots, 0, 2), (0, 0, \cdots, 0, 3, 0, \cdots, 0, 1)\rangle
 \end{array}
$

~~~~~~~where $1\le k\le n-1.$ \\

Suppose that $|F\cap (BH^{n}-\partial D_{n-1})|=k$, $|F\cap \partial D_{n-1}|=t$, then $k+t\le 2n-2$ and $t\ge 1$. Hence, there exists at least one desired $8$-cycle owing to  $2((2n-2)-k)-t\ge 2(2n-2)-(k+t)-k\ge 2n-2-k\ge 1$.
$\qed$
\section{ Edge-{bipancyclicity} of $BH_n$ under edge faults}

In this section, we consider the edge-bipancyclicity of $BH_n$ for at most $(2n-2)$ faulty edges.

  Let $e=(x, y)$ be an edge between $BH^0_{n-1}$ and $BH^1_{n-1}$ and suppose that $x', y'\in BH_n$ such that $N_{BH_n}(x)=N_{BH_n}(x')$ and $N_{BH_n}(y)=N_{BH_n}(y')$. Let $F=\{(x, y'), (x', y)\}$. From the concluding remarks of \cite{Dongqincheng2}, we have that
   there does not exists a cycle of length $4$ in $BH_n\setminus F$ that contains $e$. Thus,
 in the following,  we prove that $BH_n$ is $(2n-2)$ edge-fault-tolerant $6$-bipancyclic.

\begin{lem}\label{L5} The balanced hypercube $BH_2$ is $2$-edge-fault-tolerant 6-bipancyclic.
\end{lem}

\noindent{\bf Proof: } The proof is rather long, and we therefore provide it in Appendix A.

\begin{thm}\label{TH1}  The balanced hypercube $BH_n$ is $(2n-2)$-edge-fault-tolerant edge 6-bipancyclic for $n\ge2$.
\end{thm}

\noindent{\bf Proof: } We prove this theorem by induction  on $n$. By Lemma \ref{L5}, the theorem holds for $n=2$. Assume that it is true for $2\le k< n$. Let $F$ be any subset of $E(BH_n)$ with $|F|\le 2n-2$  and $F_i=\partial D_i\cap F$ for $0\le i\le n-1$. We get $|F|=\sum_{i=0}^{n-1}|F_i|$. Accordingly, without loss of generality, we can assume that $|F_{n-1}|\ge |F_{n-2}|\ge\cdots\ge|F_0|$. Let $F^i=F\cap E(BH^{i}_{n-1})$ for $0\le i\le 3$. We obtain  $F=F^0\cup F^1\cup F^2\cup F^3\cup F_{n-1}$ and $|F^0\cup F^1\cup F^2\cup F^3|\le 2n-4$. Let $e$ be any edge in $BH_n\setminus F$ and  $\ell$ be any even integer with $6\le \ell\le 2^{2n}$. We need to construct an $\ell$-cycle in $BH_n\setminus F$ containing $e$.

\noindent{\bf Case 1: } $e=(u, v)\notin \partial D_{n-1}$.

  Without loss of generality, we can assume that $e\in BH^0_{n-1}$.

  \noindent{\bf Subcase 1.1: } $6\le \ell\le 2^{2n-2}$\vspace{1.0ex}.

  Since $|F^0|\le |F^0\cup F^1\cup F^2\cup F^3|\le 2n-4$, by induction hypothesis, it holds\vspace{1.0ex}.

  \noindent{\bf Subcase 1.2: } $2^{2n-2}+2\le \ell\le 2^{2n-1}+6$.\vspace{1.0ex}

By induction hypothesis, there exists a  fault-free Hamiltonian cycle $C$ in $BH^0_{n-1}$ containing $e$, say $\langle c^1, c^2, \cdots$, $c^{2^{2n-2}}, c^1\rangle$ where $c^1=u, c^{2^{2n-2}}=v$. \vspace{1.0ex} We can observe that $C\setminus\{e\}$ is a $(2^{2n-2}-1)$-path. Then, $M=\{(c^1, c^{2^{2n-3}+6}), \cdots$, $(c^i, c^{2^{2n-3}+i+5}), \cdots$, $(c^{2^{2n-3}-5}$, $c^{2^{2n-2}})\}$ is a set with $2^{(2n-3)}-5$  \vspace{1.0ex} pairs of distinct vertices of $BH^0_{n-1}$ such that $d_{C\setminus\{e\}}(c^i, c^{2^{2n-3}+i+5})=2^{2n-3}+5$ for all $1\le i\le 2^{2n-3}-5$. Thus, $c^i$ and $c^{2^{2n-3}+i+5}$ are in different partite sets. There exists at least one pair $(c^t, c^{2^{2n-3}+t+5})$ in $M$ such that

$\begin{array}{l}
  |F\cap \{e_1, e_2\mid e_1, e_2 {\rm~are~ two~} (n-1) {\rm\text{-}dimensional ~edges ~incident ~with~} c^t\}|\le 1~{\rm and}\vspace{1.0ex} \\
  |F\cap \{e_3, e_4\mid e_3, e_4 {\rm~are~ two~} (n-1) {\rm\text{-}dimensional ~edges ~incident ~with~} c^{2^{2n-3}+t+5}\}|\le 1\vspace{1.0ex}
\end{array}
$\\
owing to $2\cdot(2^{2n-3}-5)>2n-2$ for all $n\ge3$. Without loss of generality, let $c^t$ be a white vertex and $c^{2^{2n-3}+t+5}$ be a black vertex.
Then, there exist two fault-free $(n-1)$-dimensional edges $(c^t, v^1), (c^{2^{2n-3}+t+5}, u^3)$ where $v^1\in BH^1_{n-1}$ and $u^3\in BH^3_{n-1}$.  Let $P_0=\langle c^{2^{2n-3}+t+5}, c^{2^{2n-3}+t+6}, \cdots$, $c^{2^{2n-2}-1}, v, u, c^2, \cdots, c^t\rangle$.
Thus, $P_0$ is a $(2^{2n-3}-5)$-path that contains $(u, v)$.
By Lemma \ref{L6},
there exists a fault-free $2$-path $\langle u^3, v^3, u^2\rangle$ and a fault-free $2$-path $\langle v^1, u^1, v^2\rangle$ where $u^i, v^i\in BH^{i}_{n-1}$ for $1\le i\le 3$. Since $|F^2|\le 2n-4$,  by induction hypothesis, there exists a Hamiltonian cycle $C_2$ in $BH^2_{n-1}\setminus F$.
Thus,  there exist two fault-free path $P'_2, P''_2$ in $BH^2_{n-1}$ joining $u^2$ and $v^2$ with length $|V(P'_2)|$ and $2^{2n-2}-|(V(P'_2))|$, respectively, where $1\le |V(P'_2)|\le 2^{2n-3}-1$.

\noindent {\bf Subcase 1.2.1: } $|V(P'_2)|=2^{2n-3}-1$ \vspace{1.0ex}.

We can represent $\ell=\ell_0+\ell_1+\ell_2+\ell_3+4$, where $\ell_i$ satisfies one of the following conditions for $i=0, 1, 2, 3$.\vspace{1.0ex}

$\begin{array}{lllll}
   \ell_0=2^{2n-3}-5,& \ell_1=1,& \ell_2=2^{2n-3}+1, &\ell_3=1~&\rm{or} \\
   \ell_0=2^{2n-3}-5, &5\le\ell_1\le 2^{2n-2}-1, &\ell_2=2^{2n-3}-1, & \ell_3=1~&\rm{or} \\
    \ell_0=2^{2n-3}-5,& 5\le \ell_1\le 2^{2n-2}-1, &\ell_2=2^{2n-3}-1, &5\le \ell_3\le 2^{2n-2}-1.
 \end{array}
$

Since $|F^i|\le 2n-4$ for $i=1, 3$,  by the induction hypothesis, there exists an $(\ell_i+1)$-cycle $C_i$ in $BH^i_{n-1}\setminus F$ containing $(u^i, v^i)$  if $5\le \ell_i\le 2^{2n-2}-1$ for $i=1, 3$. Let\vspace{1.0ex}

$\begin{array}{lll}
 P_0=\langle c^{2^{2n-3}+t+5}, c^{2^{2n-3}+t+6}, \cdots, c^{2^{2n-2}-1}, v, u, c^2, \cdots, c^t\rangle{\rm~ with~ length}~ \ell_0 ,\vspace{1.0ex} \\
     P_1=\left\{
\begin{aligned}
  & (v^1, u^1)&\rm{if} &~~\ell_1=1,\\
 &  C_1-(v^1, u^1)&\rm{if} &~~5\le\ell_1\le 2^{2n-2}-1,
\end{aligned}
\right.\vspace{1.0ex} \\
    P_2=\left\{
\begin{aligned}
  & P'_2& \rm{if}&~~\ell_2=2^{2n-3}-1,\\
 &  P''_2&\rm{if}&~~\ell_2=2^{2n-3}+1,
\end{aligned}
\right.\vspace{1.0ex} \\
     P_3=\left\{
\begin{aligned}
  & (v^3, u^3)&\rm{if}&~~\ell_3=1,\\
 &  C_3-(v^3, u^3)&\rm{if}&~~5\le \ell_3\le 2^{2n-2}-1,
\end{aligned}
\right.\vspace{1.0ex}
 \end{array}
$

\begin{figure}[!htbp]
  \centering
  \includegraphics[width=0.41\textwidth]{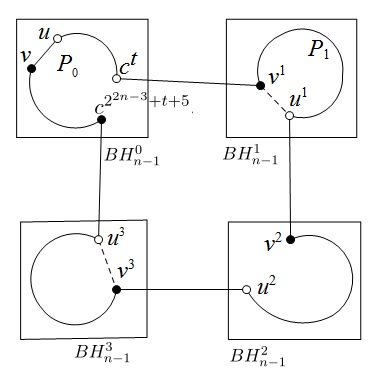}\\
  \caption{Illustration for the cycle $C$ of subcase 1.2.1 and subcase 1.2.2 in theorem \ref{TH1}. }\label{1_2_1}
\end{figure}

\noindent Then, $C=\langle c^t, v^1, P_1,  u^1, v^2, P_2, u^2, v^3, P_3, u^3, c^{2^{2n-3}+t+5}, P_0,  c^t\rangle$ ( see figure \ref{1_2_1}) is the desired cycle.

\noindent {\bf Subcase 1.2.2: } $1\le |V(P'_2)|\le 2^{2n-3}-3$ \vspace{1.0ex}.

We can represent $\ell=\ell_0+\ell_1+\ell_2+\ell_3+4$, where $\ell_i$ satisfies one of the following conditions for $i=0, 1, 2, 3$.\vspace{1.0ex}

$\begin{array}{lllll}
   \ell_0=2^{2n-3}-5, &5\le \ell_1\le 2^{2n-2}-1, &\ell_2=|V(P'_2)|, &\ell_3=1~&\rm{or} \\
    \ell_0=2^{2n-3}-5,& 5\le \ell_1\le 2^{2n-2}-1, &\ell_2=|V(P'_2)|, &5\le \ell_3\le 2^{2n-2}-1.\vspace{1.0ex}
 \end{array}
$

Since $|F^i|\le 2n-4$ for $i=1, 3$,  by the induction hypothesis, there exists an $(\ell_i+1)$-cycle $C_i$ in $BH^i_{n-1}\setminus F$ containing $(u^i, v^i)$  if $5\le \ell_i\le 2^{2n-2}-1$ for $i=1, 3$. Let\vspace{1.0ex}

$\begin{array}{lll}
 P_0=\langle c^{2^{2n-3}+t+5}, c^{2^{2n-3}+t+6}, \cdots, c^{2^{2n-2}-1}, v, u, c^2, \cdots, c^t\rangle{~\rm ~with~length~} \ell_0,\vspace{1.0ex} \\
     P_1=\left\{
\begin{aligned}
  & (v^1, u^1)&\rm{if} &~~\ell_1=1,\\
 &  C_1-(v^1, u^1)&\rm{if} &~~5\le\ell_1\le 2^{2n-2}-1,
\end{aligned}
\right.\vspace{1.0ex} \\
    P_2=P'_2, \vspace{1.0ex}\\
 P_3=\left\{
\begin{aligned}
  & (v^3, u^3)&\rm{if}&~~\ell_3=1,\\
 &  C_3-(v^3, u^3)&\rm{if}&~~5\le\ell_3\le 2^{2n-2}-1,
\end{aligned}
\right.\vspace{1.0ex}
 \end{array}
$

\noindent Then, $C=\langle c^t, v^1, P_1,  u^1, v^2, P_2, u^2, v^3, P_3, u^3, c^{2^{2n-3}+t+5}, P_0,  c^t\rangle$ ( see figure \ref{1_2_1}) is the desired cycle.

\noindent {\bf Subcase 1.3: } $2^{2n-1}+8\le \ell\le 2^{2n}$.

 We can represent $\ell=\ell_0+\ell_1+\ell_2+\ell_3+4$, where $\ell_i$ satisfies one of the following conditions for $i=0, 1, 2, 3$.\vspace{1.0ex}

$\begin{array}{lllll}
   \ell_0=2^{2n-2}-1, & 5\le\ell_1\le 2^{2n-2}-1, &\ell_2=2^{2n-2}-1, &\ell_3=1~&\rm{or} \\
    \ell_0=2^{2n-2}-1,& 5\le\ell_1\le 2^{2n-2}-1, &\ell_2=2^{2n-2}-1, &5\le \ell_3\le 2^{2n-2}-1.\vspace{1.0ex}
 \end{array}
$

By the induction hypothesis, there exists a fault-free Hamiltonian cycle $C_0$ in $BH^0_{n-1}$ containing $e$, say $\langle c^1, c^2, \cdots, c^{2^{2n-2}}, c^1\rangle$ with $c^1=u, c^{2^{2n-2}}=v$. Let $M=\{(c^1, c^2), \cdots, (c^{2i-1}, c^{2i}), \cdots$, $(c^{2^{2n-2}-1}, c^{2^{2n-2}}) \}$, then $M$ is a  set with $2^{2n-3}$ mutually disjoint edges. There exists an edge $(c^{2t-1}, c^{2t})$ in $M$ such that\vspace{1.0ex}

$ \begin{array}{l}
  |F\cap \{e_1, e_2\mid e_1, e_2 {\rm~are~ two~} (n-1) {\rm\text{-}dimensional ~edges ~incident ~with~} c^{2t-1}\}|\le 1 {\rm ~and} \vspace{1.0ex} \\
   |F\cap \{e_3, e_4\mid e_3, e_4 {\rm~are~ two~} (n-1) {\rm\text{-}dimensional ~edges ~incident ~with~} c^{2t}\}|\le 1\vspace{1.0ex}
 \end{array}$\\
since $2\cdot(2^{2n-3})>2n-2$ for all $n\ge3$.
Let $(c^{2t-1}, v^1), (c^{2t}, u^3)$ be two fault-free $(n-1)$-dimensional edges where $v^1\in BH^1_{n-1}$, $u^3\in BH^3_{n-1}$.
 By Lemma \ref{L6}, there exists a fault-free $2$-path $\langle u^3, v^3, u^2\rangle$ and a fault-free $2$-path $\langle v^1, u^1, v^2\rangle$ where $v^i, u^i\in BH^i_{n-1}$ for $i=1, 2, 3$. By Lemma \ref{LL2}, there exists a Hamiltonian path $P_2$ in $BH^2_{n-1}\setminus F$ joining $v^2$ to $u^2$. Note that $|F^i|\le 2n-4$, by the induction hypothesis, there exists an $(\ell_i+1)$-cycle $C_i$ in $BH^i_{n-1}\setminus F$ containing $(u^i, v^i)$ where $5\le\ell_i\le 2^{2n-2}-1$ for $i=1, 3$. Let \vspace{1.0ex}

$\begin{array}{lll}
 P_0=C_0-(c^{2t-1}, c^{2t}),\vspace{1.0ex} \\
     P_1=\left\{
\begin{aligned}
  & (v^1, u^1)&\rm{if} &~~\ell_1=1,\\
 &  C_1-(v^1, u^1)&\rm{if} &~~5\le\ell_1\le 2^{2n-2}-1,
\end{aligned}
\right.\vspace{1.0ex} \\
    P_2 {\rm~be~ the~ Hamiltonian~ path~ of~} BH^2_{n-1} {\rm~joining~}  v^2 {\rm~to~} u^2, \vspace{1.0ex}\\
 P_3=\left\{
\begin{aligned}
  & (v^3, u^3)&\rm{if}&~~\ell_3=1,\\
 &  C_3-(v^3, u^3)&\rm{if}&~~5\le\ell_3\le 2^{2n-2}-1,
\end{aligned}
\right.\vspace{1.0ex}
 \end{array}
$

\begin{figure}[!htbp]
  \centering
  \includegraphics[width=0.41\textwidth]{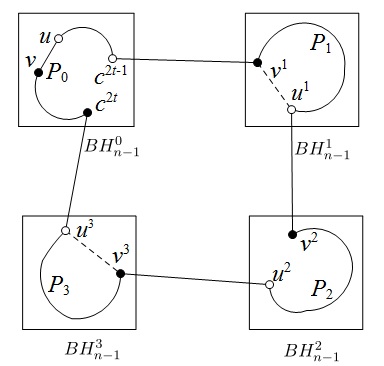}\\
  \caption{Illustration for the cycle $C$ of subcase 1.3 in theorem \ref{TH1}. }\label{1_3}
\end{figure}

\noindent Then, $C=\langle c^{2t-1}, v^1, P_1, u^1, v^2, P_2, u^2, v^3, P_3, u^3, c^{2t}, P_0, c^{2t-1}\rangle$ (see figure \ref{1_3}) is the desired cycle.

\noindent {\bf Case 2: } $e=(u, v)\in \partial D_{n-1}$\vspace{1.0ex}.

\noindent {\bf Subcase 2.1: } $|F_{n-1}|\le 2n-3$\vspace{1.0ex}.

We divide $BH_n$ into four parts, $BH^{n-2, 0}_{n-1}, BH^{n-2, 1}_{n-1}, BH^{n-2, 2}_{n-1}, and BH^{n-2, 3}_{n-1}$.
 If $|F|\le 2n-3$, then, $|F\cap (\cup_{i=0}^{3}BH^{n-2, i}_{n-1})|\le 2n-3$. If $|F|=2n-2$, note that $|F_{n-1}|\ge |F_{n-2}|\ge \cdots\ge |F_0|$, $|F_{n-1}|\le 2n-3$,  we have $|F_{n-2}|= 1$. As a result, $|F\cap (\cup_{i=0}^{3}BH^{n-2, i}_{n-1})|\le 2n-3$.

\noindent {\bf Subcase 2.1.1: } $|F\cap BH^{n-2, i}_{n-1}|\le 2n-4$ for all $i=0, 1, 2, 3$. \vspace{1.0ex}

By a similar discussion as case $1$, we obtain the result.

\noindent {\bf Subcase 2.1.2: } There exists an $i\in \{0, 1, 2, 3\}$ such that $|F\cap BH^{n-2, i}_{n-1}|=2n-3$. \vspace{1.0ex}

Without loss of generality, we can assume that $|F\cap BH^{n-2, 0}_{n-1}|=2n-3$. Thus, $|F\cap (BH_n\setminus BH^{n-2, 0}_{n-1})|\le 1$ and $|F\cap BH^{n-2, i}_{n-1}|=0$ for $i=1, 2, 3$.

\noindent {\bf Subcase 2.1.2.1: } $e\in BH^{n-2, 0}_{n-1}$. \vspace{1.0ex}

\noindent {\bf Subcase 2.1.2.1.1: } $\ell=6$. \vspace{1.0ex}

\begin{figure}[!htbp]
  \centering
  \includegraphics[width=0.45\textwidth]{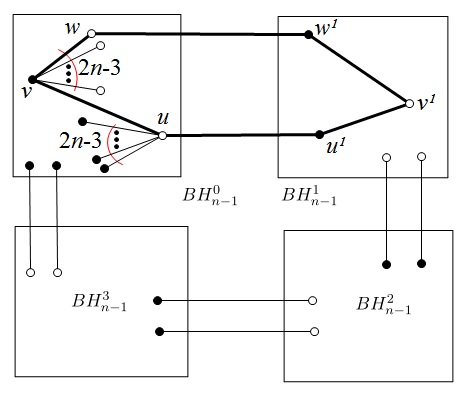}\\
  \caption{Illustration for the cycle $C$ of subcase 2.1.2.1.1 in theorem \ref{TH1}}\label{2_1_2_1_1}
\end{figure}

Note that $e$ is a fault-free edge and there are $(4n-6)$ different $2$-paths in $BH^{n-2, 0}_{n-1}$ containing $e$. Since $4n-6-(2n-3)=2n-3\ge 1$, there exists at least one fault-free $2$-path in $BH^{n-2, 0}_{n-1}$ containing $e$, say $\langle u, v, w\rangle$. Without loss of generality, let $v$ be a black vertex and $u, w$ be two white vertices. Notice that $|F_{n-2}|\le 1$, we obtain that there exists two fault-free $(n-2)$-dimensional edges $(u, u^1), (w, w^1)$ where $u^1, w^1\in BH^{n-2, 1}_{n-1}$. It is easy to verify that $d(u^1, w^1)=2$. Suppose that $v^1$ is the vertex that is adjacent to both $w^1$ and $u^1$.
 Since $|F\cap BH^{n-2, 1}_{n-1}|=0$, then $C=\langle u, v, w, w^1, v^1, u^1, u\rangle$(see figure \ref{2_1_2_1_1}) is the desired cycle.

\noindent {\bf Subcase 2.1.2.1.2: } $\ell=8$. \vspace{1.0ex}

By Lemma \ref{L4}, there are two $8$-cycles $C_1, C_2$ in $BH_n$ containing $e$ such that $E(C_1)\cap E(C_2)=e$ and $|E(C_i)\cap BH^{n-2, j}_{n-1}|=1$ for $i=1, 2, j=0, 1, 2, 3$. Note that $|F\cap (BH_n\setminus BH^{n-2, 0}_{n-1})|\le 1$. There exists at least one fault-free 8-cycle that contains $e$, say $\langle u,  v^1, u^1, v^2,  u^2, v^3, u^3, v, u\rangle$ where $u^i, v^i\in BH^{n-2, i}_{n-1}$ for $1\le i\le 3$.

\noindent {\bf Subcase 2.1.2.1.3: } $10\le \ell\le 3\cdot2^{2n-2}+2$. \vspace{1.0ex}

We can represent $\ell=\ell_0+\ell_1+\ell_2+\ell_3+4$, where $\ell_i$ satisfies one of the following conditions for $i=0, 1, 2, 3$.\vspace{1.0ex}

$\begin{array}{llll}
     \ell_0=1, 3\le\ell_1\le 2^{2n-2}-1, &\ell_2=1, &\ell_3=1~&\rm{or} \\
     \ell_0=1, 3\le\ell_1\le 2^{2n-2}-1, &3\le\ell_2\le2^{2n-2}-1, & \ell_3=1 &\rm{or}\\
     \ell_0=1, 3\le\ell_1\le 2^{2n-2}-1, &3\le\ell_2\le2^{2n-2}-1, &3\le\ell_3\le2^{2n-2}-1.\vspace{1.0ex}
 \end{array}
$

  Let $\langle u,  v^1, u^1, v^2,  u^2, v^3, u^3, v, u\rangle$ be a fault-free 8-cycle where $u^i, v^i\in BH^{n-2, i}_{n-1}$ for $1\le i\le 3$.
Since $|F\cap BH^{n-2, i}_{n-1}|=0$ for $i=1, 2, 3$. By Lemma \ref{L3}, there exists an  $(\ell_i+1)$-cycle $C_i$ in $BH^{n-2, i}_{n-1}$ containing $(u^i, v^i)$ where $3\le \ell_i\le2^{2n-2}-1$ for $1\le i\le 3$. Let\vspace{1.0ex}

$\begin{array}{lll}
     P_1=\left\{
\begin{aligned}
  & (v^1, u^1)&\rm{if} &~~\ell_1=1,\\
 &  C_1-(v^1, u^1)&\rm{if} &~~3\le\ell_1\le 2^{2n-2}-1,
\end{aligned}
\right.\vspace{1.0ex} \\
 P_2=\left\{
\begin{aligned}
  & (v^2, u^2)&\rm{if} &~~\ell_2=1,\\
 &  C_2-(v^2, u^2)&\rm{if} &~~3\le\ell_2\le 2^{2n-2}-1,
\end{aligned}
\right.\vspace{1.0ex} \\
 P_3=\left\{
\begin{aligned}
  & (v^3, u^3)&\rm{if}&~~\ell_3=1,\\
 &  C_3-(v^3, u^3)&\rm{if}&~~3\le\ell_3\le 2^{2n-2}-1,
\end{aligned}
\right.\vspace{1.0ex}
 \end{array}
$

\begin{figure}[!htbp]
  \centering
  \includegraphics[width=0.4\textwidth]{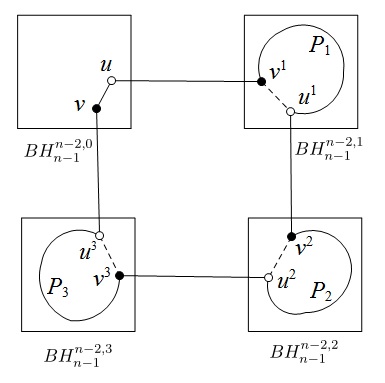}\\
  \caption{Illustration for the cycle $C$ of subcase 2.1.2.1.3 in theorem \ref{TH1}}\label{2_1_2_1_3}
\end{figure}

Then, $C=\langle u,  v^1, P_1, u^1, v^2, P_2, u^2, v^3, P_3, u^3, v, u\rangle$ (see figure \ref{2_1_2_1_3}) forms the desired cycle.

\noindent {\bf Subcase 2.1.2.1.4: } $3\cdot2^{2n-2}+4\le \ell\le2^{2n} $. \vspace{1.0ex}

We can represent $\ell=\ell_0+\ell_1+\ell_2+\ell_3+4$, where $\ell_0= 2^{2n-2}-1, 3\le \ell_i\le 2^{2n-2}-1$ for $i=1, 2, 3$.\vspace{1.0ex}

Let $\bar{e}=(u^0, v^0)$ be any faulty edge in $BH^{n-2, 0}_{n-1}$. By the induction hypothesis, there exists a Hamiltonian cycle  $C_0$ in $BH^{n-2, 0}_{n-1}-F+\{\bar{e}\}$ containing $e$. Obviously,  $|F\cap E(C_0)|\le 1$.
If $|F\cap E(C_0)|= 1$, then $\bar{e}\in E(C_0)$, we can assume that $(a^0, b^0)=\bar{e}$.
If $|F\cap E(C_0)|= 0$, let $(a^0, b^0)$ be any edge in $E(C_0)\setminus\{e\}$.
Note that $|F\cap (BH_n\setminus BH^{n-2, 0}_{n-1})|\le 1$, by Lemma \ref{L4}, there exists a fault-free $8$-cycle $\langle a^0, b^1, a^1, b^2, a^2, b^3, a^3, b^0, a^0\rangle$ in $BH_n$ where $a^i, b^i\in BH^{n-2, i}_{n-1}$ for $i=0, 1, 2, 3$.
Note that $|F\cap BH^{n-2, i}_{n-1}|=0$ for $1\le i\le 3$,  by Lemma \ref{L3}, there exists an $(\ell_i+1)$-cycle $C_i$ in $BH^{n-2, i}_{n-1}$ containing $(a^i, b^i)$ where $3\le \ell_i\le 2^{2n-2}-1$ for $1\le i\le 3$.

\begin{figure}[!htbp]
  \centering
  \includegraphics[width=0.4\textwidth]{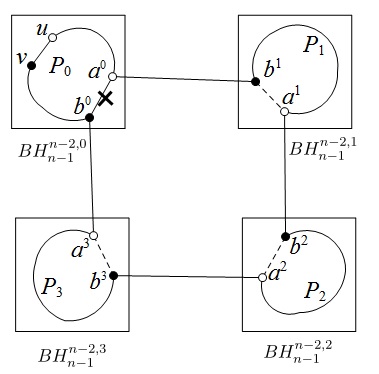}\\
  \caption{Illustration for the cycle $C$ of subcase 2.1.2.1.4 in theorem \ref{TH1}}\label{2_1_2_1_4}
\end{figure}

Let $P_i=C_i-(b^i, a^i)$ for $i=0, 1, 2, 3$, then $C=\langle  a^0, b^1, P_1, a^1, b^2, P_2, a^2, b^3, P_3, a^3, b^0, P_0, a^0\rangle$ (see figure \ref{2_1_2_1_4}) forms the desired cycle.

\noindent {\bf Subcase 2.1.2.2: } $e\in BH^{n-2, i}_{n-1}$ where $i=1, 2, 3$. \vspace{1.0ex}

Without loss of generality, we assume that $e\in BH^{n-2, 1}_{n-1}$.

\noindent {\bf Subcase 2.1.2.2.1: } $6\le \ell\le 2^{2n-2}$. \vspace{1.0ex}

Since $|F\cap BH^{n-2, 1}_{n-1}|=0$, by the induction hypothesis, it holds.

\noindent {\bf Subcase 2.1.2.2.2: } $2^{2n-2}+2\le \ell\le 2^{2n-1}-2$.\vspace{1.0ex}

We can represent $\ell=\ell_1+\ell_2+2$, where $2\le \ell_1\le 2^{2n-2}-2, \ell_2=2^{2n-2}-2$.

Since $|F\cap BH^{n-2, 1}_{n-1}|=0$, by Lemma \ref{L3}, there exists a Hamiltonian cycle $C_1$ in $BH^{n-2, 1}_{n-1}$ containing $e$, say $\langle c^0, c^1, \cdots, c^{2^{2n-2}-1}, c^0\rangle$, where $c^0=u, c^1=v$.
Let $\ell_1$ be an even integer. Then, $c^{\ell_1}$ is a white vertex and $\langle u, c^1, c^2, \cdots, c^{\ell_1}\rangle$ is an $\ell_1$-path in $BH^{n-2, 1}_{n-1}$ containing $e$ where $2\le \ell_1\le 2^{2n-2}-2$. Notice that $|F_{n-2}|\le 1$. We can assume that $(c^{\ell_1}, u^2), (u, v^2)$ are two fault-free $(n-2)$-dimensional edges where $u^2, v^2\in BH^{n-2, 2}_{n-1}$ since every vertex has two extra neighbors. By Lemma \ref{LL4}, there exists a $(2^{2n-2}-2)$-path in $BH^{n-2, 2}_{n-1}$ joining $u^2$ to $v^2$. Let \vspace{1.0ex}

\begin{figure}
  \centering
  \includegraphics[width=0.4\textwidth]{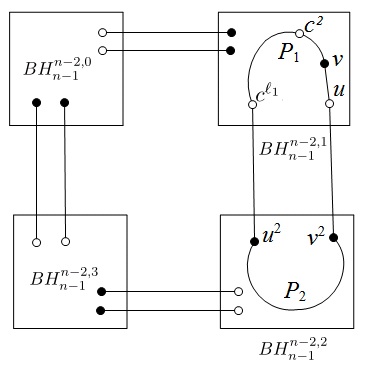}\\
  \caption{Illustration for the cycle $C$ of subcase 2.1.2.2.2 in theorem \ref{TH1}}\label{2_1_2_2_2}
\end{figure}

$\begin{array}{l}
   P_1=\langle u, v, c^2, c^3, \cdots, c^{\ell_1}\rangle, \\
   P_2~{\rm be~} {\rm the~ path~of~length~}2^{2n-2}-2~ {\rm in}~ BH^{n-2, 2}_{n-1}~{\rm joining~}u^2~{\rm and~} v^2.\vspace{1.0ex}
 \end{array}
$\\
Then, the cycle  $C=\langle u,  P_1, c^{\ell_1}, u^2, P_2, v^2, u\rangle$ (see figure \ref{2_1_2_2_2} )forms the desired cycles.

\noindent {\bf Subcase 2.1.2.2.3: } $2^{2n-1}\le \ell\le 2^{2n-1}+8$.\vspace{1.0ex}

 We can represent $\ell=\ell_0+\ell_1+\ell_2+\ell_3+4$, where $\ell_0=1, \ell_1=5, \ell_2=2^{2n-2}-1, 2^{2n-2}-9\le\ell_3\le 2^{2n-2}-1$.

Let $\langle u, v, w^1, x^1, y^1, z^1\rangle$ be a fault-free $5$-path of $BH^{n-2, 1}_{n-1}$ and
 $(z^1, u^0), (u, v^2)$ be two fault-free $(n-2)$-dimensional edges where $u^0\in BH^{n-2, 0}_{n-1}, v^2\in BH^{n-2, 2}_{n-1}$.
By Lemma \ref{L9}, there exists a $2$-path $\langle u^0, v^0, u^3\rangle$ and a $2$-path $\langle u^3, v^3, u^2\rangle$ where $u^i, v^i\in BH^{n-2, i}_{n-1}$. By Lemma \ref{L3}, there exists a $(\ell_3+1)$-cycle of $BH^{n-2, 3}_{n-1}$ containing $(u^3, v^3)$ where $2^{2n-2}-5\le\ell_3\le2^{2n-2}-1$.

\begin{figure}[!htbp]
  \centering
  \includegraphics[width=0.4\textwidth]{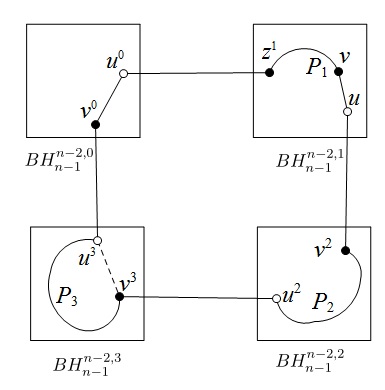}\\
  \caption{Illustration for the cycle $C$ of subcase 2.1.2.2.3 in theorem \ref{TH1}}\label{2_1_2_2_3}
\end{figure}

By Lemma \ref{LL2}, there exists a Hamiltonian path $P_2$ in $BH^{n-2, 2}_{n-1}$ joining $u^2$ and $v^2$.
Let  $P_1= \langle u, v, w^1, x^1, y^1, z^1\rangle,  P_3=C_3-(u_3, v_3)$.
Then,  $C=\langle u, P_1, z^1, u^0, v^0, u^3, P_3, v^3, u^2, P_2, v^2, u\rangle$(see figure \ref{2_1_2_2_3}) is the desired cycle.

\noindent {\bf Subcase 2.1.2.2.4: } $2^{2n-1}+10\le \ell\le 2^{2n}$.\vspace{1.0ex}

We can represent $\ell=\ell_0+\ell_1+\ell_2+\ell_3+4$ where
$
    5\le \ell_0\le 2^{2n-2}-1, \ell_1= 2^{2n-2}-1, \ell_2=2^{2n-2}-1, 3\le\ell_3\le2^{2n-2}-1.\vspace{1.0ex}
$

Let $(a^0, b^0)$ be a faulty edge in $BH^{n-2, 0}_{n-1}$, where $a^0$ is a white vertex. Assume that $(a^0, b^1), (b^0, a^3)$, $(a^3, b^3)$, and $(b^3, a^2)$ are fault-free edges where $a^i, b^i\in BH^{n-2, i}_{n-1}$ for $i=0, 1, 2, 3$.
Note that $|F\cap BH^{n-2, 1}_{n-1}|=0$, by Lemma \ref{L3},
there exits a Hamiltonian cycle $C_1$ in $BH^{n-2, 1}_{n-1}$ containing $e$.
Suppose that $N_{C_1}(b^1)=\{a^1, c^1\}$. Thus, $(b^1, a^1)\not=e$ or $(b^1, c^1)\not=e$.
Without loss of generality,  assume that $(b^1, a^1)\not=e$.
Note that $|N_{BH^{n-2, 2}_{n-1}}(a^1)|=2$ and $|F_{n-2}|\le 1$.
Suppose that $(a^1, b^2)$ is a fault-free edge where $b^2\in BH^{n-2, 2}_{n-1}$.
By Lemma \ref{LL2}, there exists a fault-free Hamiltonian path $P_2$ in $BH^{n-2, 2}_{n-1}$ joining $a^2$ and $b^2$.
By the induction hypothesis, there exists an $(\ell_0+1)$-cycle $C_0$ in $BH^{n-2, 0}_{n-1}-F+{(a^0, b^0)}$ containing $(a^0, b^0)$ where $5\le \ell_0\le 2^{2n-2}-1$.
By Lemma \ref{L3}, there exists an $(\ell_3+1)$-cycle $C_3$ in $BH^{n-2, 3}_{n-1}$ containing $a^3, b^3$ where $3\le\ell_3\le 2^{2n-2}-1$.
Let\vspace{1.0ex}

$\begin{array}{lll}
P_0= C_0-(a^0, b^0), \\
P_1=C_1-(b_1, a_1), \\
 P_2 ~{\rm be~ the~ Hamiltonian~path~joining~}a^2~{\rm and~}b^2,~\\
 P_3= C_3-(a^3, b^3).
\vspace{1.0ex}
 \end{array}
$

\begin{figure}[!htbp]
  \centering
  \includegraphics[width=0.4\textwidth]{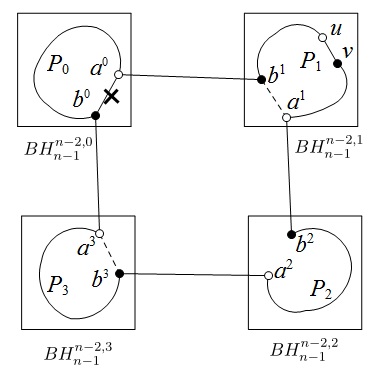}\\
  \caption{Illustration for the cycle $C$ of subcase 2.1.2.2.4 in theorem \ref{TH1}}\label{2_1_2_2_4}
\end{figure}

\noindent Then, $C=\langle a^0, P_0, b^0, a^3, P_3, b^3, a^2, P_2, b^2, a^1, P_1, b^1, a^0\rangle$ (see figure \ref{2_1_2_2_4}) is the desired cycle.

\noindent {\bf Subcase 2.2: } $|F_{n-1}|= 2n-2$ . \vspace{1.0ex}

\noindent {\bf Subcase 2.2.1: } $\ell=6$.\vspace{1.0ex}

Without loss of generality, we can assume that $e\in BH^{n-2, 0}_{n-1}$. If $|F\cap BH^{n-2, 0}_{n-1}|\le 2n-4$,  by the induction hypothesis, there exists a $6$-cycle in $BH^{n-2, 0}_{n-1}$. Thus, we assume that $|F\cap BH^{n-2, 0}_{n-1}|= 2n-3$ or $2n-2$.
Note that $e$ is a fault-free edge and there are $(4n-6)$ different $2$-paths in $BH^{n-2, 0}_{n-1}$ containing $e$. Since $4n-6-(2n-2)=2n-4\ge 1$, there exists at least one fault-free $2$-path in $BH^{n-2, 0}_{n-1}$ containing $e$, say $\langle u, v, w\rangle$. Without loss of generality, let $v$ be a black vertex and $u, w$ be two white vertices. Notice that $|F_{n-2}|=0$, we can assume that $(u, u^1), (w, w^1)$ are two  fault-free $(n-2)$-dimensional edges where $u^1, w^1\in BH^{n-2, 1}_{n-1}$. It is easy to check that $d(u^1, w^1)=2$. Suppose that $v^1$ is the vertex that is adjacent to both $w^1$ and $u^1$. Let $\bar{v}^1$ be the vertex such that $v^1$ and $\bar{v}^1$ differ in  only the inner index. Then, $(w^1, \bar{v}^1), (u^1, \bar{v}^1)\in E(BH_n)$.
 Since $|F\cap BH^{n-2, 1}_{n-1}|\le 1$, then $\langle u, v, w, w^1, v^1, u^1, u\rangle$ or $\langle u, v, w, w^1, \bar{v}^1, u^1, u\rangle$ is the desired cycle.

\noindent {\bf Subcase 2.2.2: } $\ell=8$.\vspace{1.0ex}

By Lemma \ref{L9}, there exists a fault-free $8$-cycle $\langle u=u^0, v=v^1, u^1, v^2, u^2, v^3, u^3, v^0, u\rangle$.

\noindent {\bf Subcase 2.2.3: } $10\le \ell\le 2^{2n}$.\vspace{1.0ex}

We can represent $\ell=\ell_0+\ell_1+\ell_2+\ell_3+4$, where $\ell_i$ satisfies one of the following conditions for $i=0, 1, 2, 3$.\vspace{1.0ex}

$\begin{array}{lllll}
     3\le\ell_0\le 2^{2n-2}-1, & \ell_1=1, &\ell_2=1, &\ell_3=1~&{\rm or}\\
     3\le\ell_0\le 2^{2n-2}-1, &3\le \ell_1\le 2^{2n-2}-1, &\ell_2=1, &\ell_3=1~&{\rm or}\\
     3\le\ell_0\le 2^{2n-2}-1, &3\le \ell_1\le 2^{2n-2}-1, &3\le\ell_2\le 2^{2n-2}-1, &\ell_3=1~&{\rm or}\\
     3\le\ell_0\le 2^{2n-2}-1, &3\le \ell_1\le 2^{2n-2}-1, &3\le\ell_2\le 2^{2n-2}-1, &3\le\ell_3\le2^{2n-2}-1.\vspace{1.0ex}
 \end{array}
$

Note that $F_{n-1}=2n-2$, we have  $F\cap BH^{i}_{n-1}=0$ for all $i=0, 1, 2, 3$.
 By Lemma \ref{L3}, there exists an $(\ell_i+1)$-cycle $C_i$ in $BH^i_{n-1}$ containing $(u^i, v^i)$ where $3\le \ell_i\le 2^{2n-2}-1$ for $i=0, 1, 2, 3$.

\noindent Let\vspace{1.0ex}

\begin{figure}
  \centering
  \includegraphics[width=0.4\textwidth]{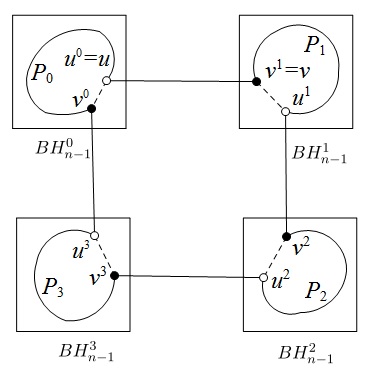}\\
  \caption{Illustration for the cycle $C$ of subcase 2.2.3 in theorem \ref{TH1}}\label{2_2_3}
\end{figure}

$\begin{array}{lll}
    P_0=\left\{
\begin{aligned}
  & (v^0, u^0)&\rm{if} &~~\ell_0=1,\\
 &  C_0-(v^0, u^0)&\rm{if} &~~3\le\ell_0\le 2^{2n-2}-1,
\end{aligned}
\right.\vspace{1.0ex} \\
     P_1=\left\{
\begin{aligned}
  & (v^1, u^1)&\rm{if} &~~\ell_1=1,\\
 &  C_1-(v^1, u^1)&\rm{if} &~~3\le\ell_1\le 2^{2n-2}-1,
\end{aligned}
\right.\vspace{1.0ex} \\
 P_2=\left\{
\begin{aligned}
  & (v^2, u^2)&\rm{if} &~~\ell_2=1,\\
 &  C_2-(v^2, u^2)&\rm{if} &~~3\le\ell_2\le 2^{2n-2}-1,
\end{aligned}
\right.\vspace{1.0ex} \\
 P_3=\left\{
\begin{aligned}
  & (v^3, u^3)&\rm{if}&~~\ell_3=1,\\
 &  C_3-(v^3, u^3)&\rm{if}&~~3\le\ell_3\le 2^{2n-2}-1,
\end{aligned}
\right.\vspace{1.0ex}
 \end{array}
$\\
Then, $C=\langle v^0, P_0, u^0, v^1, P_1, u^1, v^2, P_2, u^2, v^3, P_3, u^3, v^0\rangle~$(see figure \ref{2_2_3}) forms the desired cycle. $\qed$

\section*{Appendix A. Proof of Lemma \ref{L5}}

{\bf Lemma \ref{L5}} The balanced hypercube $BH_2$ is $2$-edge-fault-tolerant 6-bipancyclic. \vspace{1.5ex}

By Lemma \ref{LL2}, for an arbitrary fault-free edge $(u, v)$, there exists a fault-free Hamiltonian path $P$  that joins $u$ and $v$, then $\langle u, P, v, u\rangle$ is the fault-free $16$-cycle. Hence, we only need to construct a fault-free $\ell$-cycle in $BH_2$ containing $(u, v)$ where $6\le \ell\le 14$.
Suppose that $|F|=2$, without loss of generality, we can assume that $|F\cap\partial D_1|\ge|F\cap \partial D_0|$.

\noindent{\bf Case 1: } $|F\cap \partial D_1|=2, |F\cap \partial D_0|=0$. \vspace{1.0ex}

\noindent{\bf Subcase 1.1: } $e=(u, v)\in \partial D_0$. \vspace{1.0ex}

Without loss of generality, we can assume that $e=(u, v)\in BH^{0}_1$. Suppose that $u=(a_0, 0)$ is a white vertex, $v=(b_0, 0)$ is a black vertex.

\noindent{\bf Subcase 1.1.1: } $6\le \ell\le 14$. \vspace{1.0ex}

  There are three $\ell$-cycles $C_1, C_2, C_3$ in $BH_n$ containing $e$ where $6\le \ell\le 14$, such that $(C_i\cap \partial D_1)\cap(C_j\cap \partial D_1)=\emptyset$ for all $1\le i\not=j\le 3$. We list them as follows: \vspace{1.0ex}

 Three  6-cycles:\\
 $~~\left\{
\begin{aligned}
  &\langle(a_0, 0), (b_0, 0), (a_0,  3), (b_0+2,  3), (a_0+2,  3), (b_0+2, 0), (a_0, 0)\rangle; \\
 & \langle(a_0, 0), (b_0, 0), (a_0+2,  3), (b_0,  3), (a_0,  3), (b_0+2, 0), (a_0, 0)\rangle;\\
 & \langle(a_0, 0), (b_0, 0), (a_0+2, 0), (b_0,  1), (a_0,  1), (b_0+2,  1), (a_0, 0)\rangle.
\end{aligned}
\right.$\vspace{1.0ex}

 Three  8-cycles:\\
 $~~\left\{
\begin{aligned}
  &\langle(a_0, 0), (b_0, 0), (a_0,  3), (b_0,  3), (a_0,  2), (b_0+2,  3), (a_0+2,  3), (b_0+2, 0), (a_0, 0)\rangle; \\
 & \langle(a_0, 0), (b_0, 0), (a_0+2,  3), (b_0,  3), (a_0+2,  2), (b_0+2,  3), (a_0,  3), (b_0+2, 0), (a_0, 0)\rangle;\\
 & \langle(a_0, 0), (b_0, 0), (a_0+2, 0), (b_0,  1), (a_0+2,  1), (b_0,  2), (a_0,  1), (b_0+2,  1), (a_0, 0)\rangle\vspace{1.0ex}.
\end{aligned}
\right.$\vspace{1.0ex}

Three  10-cycles:\\
 $~~\left\{
\begin{aligned}
  &\langle(a_0, 0), (b_0, 0), (a_0,  3), (b_0,  3), (a_0,  2), (b_0,  2), (a_0+2,  2), (b_0+2,  3), (a_0+2,  3), (b_0+2, 0), (a_0, 0)\rangle; \\
 & \langle(a_0, 0), (b_0, 0), (a_0+2,  3), (b_0,  3), (a_0+2,  2), (b_0,  2), (a_0,  2), (b_0+2,  3),(a_0,  3), (b_0+2, 0), (a_0, 0)\rangle;\\
 & \langle(a_0, 0), (b_0, 0), (a_0+2, 0), (b_0,  1), (a_0+2,  1), (b_0+2,  2), (a_0+2,  2), (b_0,  2),  (a_0,  1), (b_0+2,  1), (a_0, 0)\rangle\vspace{1.0ex}.
\end{aligned}
\right.$\vspace{1.0ex}

Three  12-cycles:\\
 $~~\left\{
\begin{aligned}
  &\langle(a_0, 0), (b_0, 0), (a_0,  3), (b_0,  3), (a_0+2,  3), (b_0+2,  3), (a_0+2,  2), (b_0,  2), (a_0,  2), (b_0+2,  2),  (a_0,  1),\\ & (b_0,  1), (a_0, 0)\rangle; \\
 & \langle(a_0, 0), (b_0, 0), (a_0+2,  3), (b_0+2,  3), (a_0,  3), (b_0,  3), (a_0,  2), (b_0+2,  2), (a_0+2, 2), (b_0,  2) , (a_0+2,  1),\\ & (b_0+2,  1), (a_0, 0)\rangle;\\
 & \langle(a_0, 0), (b_0, 0), (a_0+2, 0), (b_0,  1), (a_0+2,  1), (b_0+2,  2), (a_0+2,  2), (b_0,  2),  (a_0,  2), (b_0+2,  3), (a_0,  3), \\ & (b_0+2, 0), (a_0, 0)\rangle\vspace{1.0ex}.
\end{aligned}
\right.$\vspace{1.0ex}

Three  14-cycles:\\
 $~~\left\{
\begin{aligned}
  &\langle(a_0, 0), (b_0, 0), (a_0,  3), (b_0,  3), (a_0+2,  3), (b_0+2,  3), (a_0+2,  2), (b_0,  2),  (a_0,  2), (b_0+2,  2), (a_0+2,  1),\\ & (b_0,  1), (a_0,  1), (b_0+2,  1), (a_0, 0)\rangle; \\
 & \langle(a_0, 0), (b_0, 0), (a_0+2,  3), (b_0+2,  3), (a_0,  3), (b_0,  3), (a_0+2,  2), (b_0+2,  2),  (a_0, 2), (b_0,  2), (a_0+2,  1), \\ & (b_0+2,  1), (a_0,  1), (b_0,  1), (a_0, 0)\rangle;\\
 & \langle(a_0, 0), (b_0, 0), (a_0+2, 0), (b_0,  1), (a_0,  1), (b_0+2,  2), (a_0+2,  2), (b_0,  2), (a_0,  2), (b_0+2,  3), (a_0,  3),\\ & (b_0,  3), (a_0+2,  3), (b_0+2, 0), (a_0, 0)\rangle\vspace{1.0ex}.
\end{aligned}
\right.$\vspace{1.0ex}

Notice that $|F\cap \partial D_1|=2, |F\cap \partial D_0|=0$,  there exists  at least one fault-free $\ell$-cycle in $BH_2$ containing $e$ where $6\le \ell\le 14$.

\noindent{\bf Subcase 1.2: } $e=(u, v)\in \partial D_1$. \vspace{1.0ex}

Without loss of generality, we can assume that $e=(u, v)=(u^0, v^1)$ is an edge between $BH^0_1$ and $BH^1_1$ where $u^0=(a_0, 0), v^1=(b_0, 1)$.

\noindent{\bf Subcase 1.2.1: } $\ell=6, 8$. \vspace{1.0ex}

There exist three $\ell$-cycles $C_1, C_2, C_3$ in $BH_2$ containing $e$  where $\ell=6$ or 8, such that $(C_i\cap \partial D_1)\cap (C_j\cap \partial D_1)=\{e\}$ for $1\le i\not=j\le 3$. We list them as follows:
\vspace{1.0ex}

 Three  6-cycles:\\
 $~~~~~~\left\{
\begin{aligned}
  &\langle(a_0, 0), (b_0,  1), (a_0,  1), (b_0+2,  1), (a_0+2, 0), (b_0, 0), (a_0, 0)\rangle; \\
 & \langle(a_0, 0), (b_0,  1), (a_0,  1), (b_0+2,  2), (a_0+2,  1), (b_0+2,  1), (a_0, 0)\rangle;\\
 & \langle(a_0, 0), (b_0,  1), (a_0+2, 0), (b_0+2, 0), (a_0,  3), (b_0, 0), (a_0, 0)\rangle.
\end{aligned}
\right.$\vspace{1.0ex}

 Three  8-cycles:\\
 $~~~~~~\left\{
\begin{aligned}
  &\langle(a_0, 0), (b_0,  1), (a_0,  1), (b_0,  2), (a_0, 2), (b_0, 3), (a_0, 3), (b_0, 0), (a_0, 0)\rangle; \\
 & \langle(a_0, 0), (b_0,  1), (a_0+2,  1), (b_0+2,  2), (a_0+2,  2), (b_0+2,  3), (a_0+2, 3), (b_0+2, 0), (a_0, 0)\rangle;\\
 & \langle(a_0, 0), (b_0,  1), (a_0,  1), (b_0+2,  2), (a_0, 2), (b_0+2, 3), (a_0, 3), (b_0+2, 0), (a_0, 0)\rangle.
\end{aligned}
\right.$\vspace{1.0ex}

Since $|F\cap \partial D_1|=2, |F\cap \partial D_0|=0$,  and $e$ is a fault-free edge, then there exists  at least one fault-free 6-cycle  and  one fault-free 8-cycle in $BH_2$ containing $e$.

\noindent{\bf Subcase 1.2.2: } $10\le\ell\le 14$. \vspace{1.0ex}

By the proof of subcase 1.2.1, there exists a fault-free $8$-cycle  $C$ that contains $e$ such that $|C\cap BH^{i}_1|=1$ for $0\le i\le 3$, say $\langle u^0, v^0, u^3, v^3, u^2, v^2, u^1, v^1, u^0\rangle$ where $u^i, v^i\in BH^{i}_{n-1}$ for $i=0, 1, 2, 3$. Since $|F\cap \partial D_0|=0$. It is easy to check that there exists an $\ell_i$-path $P_i$ in $BH^i_1$  joining $u^i$ to $v^i$ where $\ell_i=1$ or $3$ for $i=1, 2, 3$. Then, the cycle $\langle u^0, v^0, u^3, P_3, v^3, u^2, P_2, v^2, u^1, P_1, v^1, u^0\rangle$ with length $\ell=5+\ell_1+\ell_2+\ell_3$ forms the desired cycle.

\noindent{\bf Case 2: } $|F\cap \partial D_1|=1, |F\cap \partial D_0|=1$. \vspace{1.0ex}

\noindent{\bf Subcase 2.1: } $e=(u, v)\in \partial D_1$. \vspace{1.0ex}

Without loss of generality, we can assume that $e=(u, v)=(u^0, v^1)$ is an edge between $BH^0_1$ and $BH^1_1$ where $u^0=(a_0, 0)$, $v=(b_0,  1)$.

\noindent{\bf Subcase 2.1.1: } $\ell=6$. \vspace{1.0ex}

If $((a_0+2, 0), (b_0+2, 1))$ is a fault-free edge. Let
\begin{center}
$\begin{array}{l}
  C_1=\langle u^0, v^1, (a_0, 1), (b_0+2, 1), (a_0+2, 0), (b_0, 0), u^0\rangle;\vspace{1.0ex} \\
  C_2=\langle u^0, v^1, (a_0+2, 1), (b_0+2, 1), (a_0+2, 0), (b_0+2, 0), u^0\rangle.
\end{array}$
\end{center}
Then, $C_1, C_2$ are two cycles in $BH_2$ containing $e$ and $C_1\cap C_2=\{e, ((a_0+2, 0), (b_0+2, 1))\}$ is the fault-free edge set. Thus, $C^1$ or $C^2$ is a fault-free 6-cycle.

If $((a_0+2, 0), (b_0+2, 1))$ is a faulty edge. Then, $(u^0, (b_0+2, 1))$ is a fault-free edge.  Let
\begin{center}
$\begin{array}{l}
  C_3=\langle u^0, v^1, (a_0, 1), (b_0, 2), (a_0+2, 1), (b_0+2, 1), u^0\rangle;\vspace{1.0ex} \\
  C_4=\langle u^0, v^1, (a_0+2, 1), (b_0+2, 2), (a_0, 1), (b_0+2, 1), u^0\rangle.
\end{array}$
\end{center}
Then, $C_3, C_4$ are two cycles in $BH_2$ containing $e$ and $C_3\cap C_4=\{e, (u^0, (b_0+2, 1))\}$ is the fault-free edge set. Thus, $C^3$ or $C^4$ is a fault-free 6-cycle.

\noindent{\bf Subcase 2.1.2: } $\ell=8$. \vspace{1.0ex}

By Lemma \ref{L9}, it holds.

\noindent{\bf Subcase 2.1.3: } $10\le\ell\le 14$. \vspace{1.0ex}

By the proof of subcase 2.1.2, there exists a fault-free $8$-cycle $C$ that contains $e$ such that $|C\cap BH^{i}_1|=1$ for $0\le i\le 3$, say $\langle u^0, v^0, u^3, v^3, u^2, v^2, u^1, v^1, u^0\rangle$ where $u^i, v^i\in BH^{i}_{n-1}$ for $i=0, 1, 2, 3$. Note that $|F\cap \partial D_0|=1$. Without loss of generality, let $|F\cap BH^0_{1}|=1$.  It is  easy to check that there exists an $\ell_i$-path in $BH^i_1$  joining $u^i$ to $v^i$ where $\ell_i=1$ or $3$ for $i=1, 2, 3$. Then,  $\langle u, v^0, u^3, P_3, v^3, u^2, P_2, v^2, u^1, P_1, v, u\rangle$ with length $\ell=5+\ell_1+\ell_2+\ell_3$ forms the desired cycle.

\noindent{\bf Subcase 2.2: } $e=(u, v)\in \partial D_0$. \vspace{1.0ex}

We divide $BH_2$ into four  $BH_1$s, denoted by $\overline{BH^0_1}, \overline{BH^1_1}, \overline{BH^2_1}, and \overline{BH^3_1}$ , by deleting all  1-dimensional edges. Then, $e$ is an edge between $\overline{BH^i_1}$ and $\overline{BH^{i+1}_1}$ for $0\le i\le 3$. By a similar discussion for subcase 2.1, we obtain the result.
$\qed$

\section{Conclusion}

In this paper, we consider the edge-bipancyclicity of $BH_n$ for at most $(2n-2)$ faulty edges and prove that each fault-free edge lies on a fault-free cycle of any even length from $6$ to $2^{2n}$. Our result improves the results of  Hao et al. \cite{R.X.Hao} and Cheng  et al. \cite{Dongqincheng2} and it is optimal with respect to the maximum number of tolerated edge faults. In addition, it is of interest to consider the problem of fault-tolerant embedding cycles with each vertex incident to at least two non-faulty edges.

\section*{Acknowledgement}

The authors would like to express their gratitude to the anonymous referees for their kind suggestions and useful comments on the original manuscript, which resulted in this final version.
This research is supported by the National Natural Science Foundation of China (11571044, 61373021, 11461004),  the Fundamental Research Funds for the Central University of China .

\end{CJK*}

\end{document}